\documentclass[12pt]{article}
\usepackage{mathrsfs}
\usepackage{amssymb} \textwidth 150mm \textheight 220mm \oddsidemargin
15pt \evensidemargin 0pt \topmargin 0cm \headsep 0.3cm

\usepackage{amsmath}
\usepackage{graphicx}
\usepackage{indentfirst}
\usepackage{cite}
\usepackage{float}
\usepackage{CJK}
\usepackage{amsmath}
\usepackage{graphicx,cite}
\usepackage{multirow}
\usepackage{makecell}

\usepackage{subfigure}

\usepackage[justification=centering]{caption}

\begin{document}
\centerline {\Large\bf{Limit cycles appearing from perturbations}}
\vskip 0.2 true cm
\centerline {\Large\bf{of cubic piecewise smooth center}}
\centerline {\Large\bf{with double invariant real straight lines}}

\vskip 0.2 true cm

\centerline{\bf  Jihua Yang$^{a,b}$, Liqin Zhao$^{a,*}$}
\centerline{ $^a$ School of Mathematical Sciences, Beijing Normal University,}
\centerline{Laboratory of Mathematics and Complex Systems, Ministry of Education,}
\centerline{Beijing 100875, The People's Republic of China}
\centerline{ $^b$ School of Mathematics and Computer Science,}
\centerline{Ningxia Normal University, Guyuan 756000, The People's Republic of China}

\footnotetext[1]{This work was supported by the National Natural Science Foundation of China (11671040,11601250), the Fundamental Research Funds for the Central Universities, the Science and Technology Pillar Program of Ningxia(KJ[2015]26(4)) and the Key Program of Higher Education of Henan(16A110038,17B110003).\\
*Author for corresponding. E-mail: zhaoliqin@bnu.edu.cn(L. Zhao), yangjh@mail.bnu.edu.cn(J. Yang).}

\vskip 0.3 true cm

\noindent
{\bf Abstract}\, This paper investigates the exact number of limit cycles given by the averaging theory of first order for the piecewise smooth integrable non-Hamiltonian system
\begin{eqnarray*}
(\dot{x},\ \dot{y})=\begin{cases}
(-y(x+a)^2+\varepsilon f^+(x,y),\ x(x+a)^2+\varepsilon g^+(x,y)),\ \ x\geq0,\\
(-y(x+b)^2+\varepsilon f^-(x,y),\ x(x+b)^2+\varepsilon g^-(x,y)),\ ~ \, x<0,\\
\end{cases}\end{eqnarray*}
where $ab\neq 0$, $0<|\varepsilon|\ll 1$, and $f^\pm(x,y)$ and $g^\pm(x,y)$ are polynomials of degree $n$. It is proved that the exact number of limit cycles emerging from the period annulus surrounding the origin is linear depending on $n$ and it is at least twice the associated estimation of smooth systems.

\section{Introduction and the main results}

\vskip 0.2 true cm
\setcounter{equation}{0}
\renewcommand\theequation{1.\arabic{equation}}

In recent years the studying of discontinuous differential systems has become one of the  frontiers between mathematics, physics and engineering. This interest has been stimulated by discontinuous phenomena in mechanics, electrical engineering, biology, and the theory of automatic control(see, for instance, [1],[9],[16], [27] and the references therein).

\vskip 0.1 true cm

As in the smooth differential systems, one of the main problem is to determine the number of limit cycles and their distributions. Therefore, the studying of periodic solutions of discontinuous differential systems can be seen as an extension of the Hilbert's 16th problem. For the study of Hilbert's 16th problem, we refer [4],[12],[15],[17],[18] and the references therein.

\vskip 0.1 true cm

The piecewise differential systems with two zones separated by a straight line is the simplest case, and they can be expressed as follows:
\begin{eqnarray}
(\dot{x},\ \dot{y})=\begin{cases}
(p^+(x,y),\ q^+(x,y)),\ \ x>0,\\
(p^-(x,y),\ q^-(x,y)),\ \ x<0,\\\end{cases}\end{eqnarray}
where $p^\pm$ and $q^\pm$ are $C^\infty$ functions. We call system (1.1) to be  a linear(resp.quadratic or cubic) piecewise smooth differential system if each of its sub-systems is linear(resp.quadratic or cubic).

\vskip 0.2 true cm

For the general piecewise smooth differential systems (1.1), Coll {\it et al.} [5] obtained the expressions for the first three Lyapunov constants.
In [7], [8] and [14], the authors investigated the number of limit cycles for  linear piecewise smooth differential systems.
In [3], [23] and [25] the authors studied the number of limit cycles or the bounded solutions for quadratic piecewise smooth differential systems.
In [13], [20-21] and [28-29], the authors discussed the number of limit cycles bifurcating from the the period annulus, or near the origin, or near the generalized polycycles of the  piecewise smooth Hamiltonian systems.

 \vskip 0.2 true cm
In [23], Llibre and Mereu studied the maximum number of limit cycles given by the averaging theory of first order for discontinuous differential systems, which can bifurcate from the periodic orbits of the quadratic isochronous centers
\begin{eqnarray*}\dot {x}=-y+x^2,~~\dot {y}=x+xy\end{eqnarray*}
and
\begin{eqnarray}\dot {x}=-y+x^2-y^2,~~\dot {y}=x+2xy\end{eqnarray}
when they are perturbed insider the class of
all discontinuous quadratic polynomial differential systems with the straight line
of discontinuity $y=0$.

\vskip 0.2 true cm

In [19], Li and Liu studied the following piecewise smooth perturbed integrable differential system
\begin{eqnarray}
(\dot{x},\ \dot{y})=
\begin{cases}
(-y(a x+1)+\varepsilon f^+(x,y),\ x(a x+1))+\varepsilon g^+(x,y),\ \ x\geq0,\\
(-y(b x+1)+\varepsilon f^-(x,y),\ x(b x+1))+\varepsilon g^-(x,y),\ ~ \, x<0,\\
\end{cases}\end{eqnarray}
where $a, b\in R$, $0<|\varepsilon|\ll 1$, and $f^{\pm}$ and $g^{\pm}$ are arbitrary polynomials of degree $n$.

\vskip 0.2 true cm

In [10], using the averaging method, Gin\'{e} and Llibre studied the maximum number of limit cycles that can bifurcate from the period annulus surrounding the origin of the
following cubic system
\begin{eqnarray} \dot{x}=-yh(x,y), ~~\dot{y}=xh(x,y)\end{eqnarray}
under the perturbations of arbitrary cubic polynomials, where
$h(x,y)=0$ is a conic such that $h(0,0)\neq0$. The conics in $\mathbb{R}^2$ are classified as ellipses, complex ellipses, hyperbolas, two complex straight lines intersecting in real point, two intersecting real straight lines, parabolas, two real parallel straight lines, two complex parallel straight lines and one double real straight line. More precisely, a conic $h=h(x,y)=0$ is given by one of the following nine cases:
\vskip 0.2 true cm

(E) Ellipse $h=(x+\alpha)^2+(y+\beta)^2-1=0$ with $\alpha^2+\beta^2\neq1$.

(CE) Complex ellipse $h=(x+\alpha)^2+(y+\beta)^2+1=0$.

(H) Hyperbola $h=(x+\alpha)^2-y^2-1=0$ with $\alpha^2\neq1$.

(CL) Two complex straight lines intersecting in a real point $h=(x+\alpha)^2+(y+\beta)^2=0$ with $\alpha\beta\neq0$.

(RL) Two real straight lines intersecting in a point (Lotka-Volterra systems) $h=(x+\alpha)(y+\beta)=0$ with $\alpha\beta\neq0$.

(P) Parabola $h=x-\alpha-y^2$ with $\alpha\neq0$.

(RPL) Two real parallel straight lines $h=(x+\alpha)^2-1=0$ with $\alpha^2\neq1$.

(CPL) Two complex parallel straight lines $h=(x+\alpha)^2+1=0$.

(DL) One double invariant real straight line $h=(x+\alpha)^2=0$ with $\alpha\neq0$.

\vskip 0.2 true cm

Motivated by [23], [19] and [10], in this paper, we intend to study
the number of limit cycles that can bifurcate from the period annulus surrounding the origin of the following cubic system
\begin{eqnarray}
(\dot{x},\ \dot{y})=
\begin{cases}
(-y(x+a)^2,~~ x(x+a)^2,\ \ x\geq0,\\
(-y(x+b)^2,~~ x(x+b)^2,\ ~ \, x<0,\\
\end{cases}
\end{eqnarray}
under the perturbations of all discontinuous  polynomial of degree $n$
 with the straight line of discontinuity $y=0$. More precisely, we will consider
the number of limit cycles given by the averaging theory of first order for the following piecewise smooth integrable differential system
\begin{eqnarray}
(\dot{x},\ \dot{y})=\begin{cases}
(-y(x+a)^2+\varepsilon f^+(x,y),\ x(x+a)^2+\varepsilon g^+(x,y)),\ \ x\geq0,\\
(-y(x+b)^2+\varepsilon f^-(x,y),\ x(x+b)^2+\varepsilon g^-(x,y)),\ ~ \, x<0,\\
\end{cases}
\end{eqnarray}
where $ab\neq0$, $0<|\varepsilon|\ll1$, and
\begin{eqnarray}
        f^\pm(x,y)=\sum\limits_{i+j=0}^na^\pm_{i,j}x^iy^j,\ \
        g^\pm(x,y)=\sum\limits_{i+j=0}^nb^\pm_{i,j}x^iy^j.\end{eqnarray}

System (1.5) has the first integral $H^\pm(x,y)=x^2+y^2$ with respect to $x\geq0$ and $x<0$, and it is a piecewise smooth integrable non-Hamiltonian system with a generalized center $(0,0)$. It has one double invariant straight line $(x+a)^2=0$ (resp. $(x+b)^2=0$) for the case $a<0$ (resp. $b>0$). Let
\begin{eqnarray}r_1:=\begin{cases}
-a,\ \ ~~~a<0,\\
+\infty,\, ~~ ~a>0,\\
\end{cases}
r_2:=\begin{cases}
b,\quad \quad~  b>0,\\
+\infty,\, ~~ ~b<0;\\
\end{cases}~~r_0:=\min\{r_1,r_2\}, \end{eqnarray}
and denote by $H(n)$ the maximum number of limit cycles of system (1.6) bifurcating  from the period annulus $\bigcup\limits_{0<h<r_0}L_h$ for all possible $f^\pm(x,y)$ and $g^\pm(x,y)$ satisfying (1.7) up to the first order averaging method. Our main result is the following theorem.

 \vskip 0.2 true cm
\noindent
{\bf Theorem 1.1.}\, {\it Consider system (1.6). Suppose that $ab\neq0$.  Then, for any $n\geq1$, we have}
\begin{eqnarray}
H(n)=\begin{cases}
4[\frac{n+1}{2}]+\frac 32(1+(-1)^n),\quad \text{if}\ a\neq-b;\\
~~~~~3[\frac{n+1}{2}]+(-1)^n, \quad\ \quad ~\text{if}\ a=-b.\end{cases}\end{eqnarray}

\vskip 0.2 true cm
\noindent
{\bf Theorem 1.2.}\, {\it System (1.4) with $h(x,y)=(x+a)^2(a\neq 0)$ has exactly
$n$ limit cycles under the perturbations of arbitrary polynomials with degree $n$ up to
first order averaging function in $\varepsilon$ bifurcating from the period annuli surrounding the origin.}

\vskip 0.2 true cm
\noindent
{\bf Remark 1.1.} (i) It is easy to see that system (1.5) corresponds to the case of (DL) defined in [10] if $a=b$.

\vskip 0.2 true cm

(ii) It was proved in [10] that system (1.4) with $h(x,y)=(x+a)^2(a\neq 0)$ has
at most 4 limit cycles under the perturbations of arbitrary cubic polynomials up to
first order averaging function in $\varepsilon$.

\vskip 0.2 true cm

(iii) Suppose that $a=b$ together with $a^+_{ij}=a^-_{ij}$ and $b^+_{ij}=b^-_{ij}$ in (1.7), then system (1.6) is smooth for all $(x,y)\in R^2$. By Theorems 1.1 and 1.2, we know that the piecewise smooth differential systems (1.6) can bifurcate at least twice the number of limit cycles than the corresponding smooth systems.

\vskip 0.2 true cm

This paper is organized as follows. In section 2, we will give some preliminaries.
In section 3, we will obtain the expression of averaged function and the lower bound
of limit cycles for system (1.6).  In section 4, we will obtain  the upper bound
of limit cycles for system (1.6) by Argument Principle. In section 5, we will prove Theorem 1.2.

\vskip 0.2 true cm

\section{Preliminaries}

We first summarize some results on the first order averaging method for discontinuous differential systems. For a general introduction to the averaging method, we recommend the readers the references [2],[22],[24] and [26].
\vskip 0.2 true cm
\noindent
{\bf Lemma 2.1.}\, {\it ([22]) Consider the following discontinuous differential equation
$$\frac{dr}{d\theta}=\varepsilon F^0(\theta,r)+\varepsilon^2 R^0(\theta,r,\varepsilon),\eqno(2.1)$$
with
\begin{eqnarray*}
\begin{aligned}
F^0(\theta,r)&=F_1(\theta,r)+\textup{sign}(\mu(\theta,r))F_2(\theta,r),\\
R^0(\theta,r,\varepsilon)&=R_1(\theta,r,\varepsilon)+\textup{sign}(\mu(\theta,r))R_2(\theta,r,\varepsilon),
\end{aligned}\end{eqnarray*}
where $F_1,F_2:\mathbb{R}\times D\rightarrow\mathbb{R}$, $R_1,R_2:\mathbb{R}\times D\times(-\varepsilon_0,\varepsilon_0)\rightarrow \mathbb{R}$, $\mu:\mathbb{R}\times D\rightarrow\mathbb{R}$ are continuous functions, T-periodic in the variable $\theta$ and $D$ is an open interval of $\mathbb{R}$. Let $\textup{sign}(\cdot)$ be the sign function.
We also suppose that $\mu$ is a $C^1$ function having 0 as a regular value. Denote by $\mathcal{M}=\mu^{-1}(0)$, by $\Sigma=\{0\}\times D\subset \mathcal{M}$, by $\Sigma_0=\Sigma\setminus \mathcal{M}\neq\varnothing$, and its elements by $z:=(0,z)\notin\mathcal{M}$.

\vskip 0.2 true cm

Define the averaged function $f^0:D\rightarrow \mathbb{R}$ as
$$f^0(r)=\int_0^TF^0(\theta,r)d\theta.\eqno(2.2)$$
We assume the following three conditions:}

(i)\, {\it $F_1$, $F_2$, $R_1$, $R_2$ and $\mu$ are locally Lipschitz with respect to $r$.}

(ii)\, {\it For $a\in\Sigma_0$ with $f^0(a)=0$, there exist a neighborhood $V$ of a such that $f^0(z)\neq0$ for all $z\in\overline{V}\setminus\{a\}$ and the Brouwer degree function $d_B(f^0,V,0)\neq0$.}

(iii)\, {\it If $\partial \mu/\partial\theta\neq0$, then for all $(\theta,r)\in\mathcal{M}$ we have that $(\partial \mu/\partial\theta)(\theta,r)\neq0$; and if $\partial \mu/\partial\theta\equiv0$, then $\langle\nabla_r\mu,F_1\rangle^2-\langle\nabla_r\mu,F_2\rangle^2>0$ for all $(\theta,z)\in[0,T]\times\mathcal{M}$, where $\nabla_r\mu$ denotes the gradient of the function $\mu$ restricted to variable $r$.}

{\it Then for $|\varepsilon|>0$ sufficiently small, there exists a T-periodic solution $r(\theta,\varepsilon)$ of system (2.1) such that $r(0,\varepsilon)\rightarrow0$ (in the sense of Hausdorff distance) as $\varepsilon\rightarrow0$.}

\vskip 0.2 true cm

In order to verify the hypothesis (ii) of Lemma 2.1, we have the following remark, see for instance [2].

\vskip 0.2 true cm
\noindent
{\bf Remark 2.1.}\, Let $f^0:D\rightarrow\mathbb{R}$ be a $C^1$ function with $f^0(a)=0$, where $D$ is an open interval of $\mathbb{R}$ and $a\in D$. Whenever the Jacobian of $J_{f^0}(a)\neq0$, there exists a neighborhood $V$ of $a$ such that $f^0(r)\neq0$ for all $r\in\overline{V}\setminus\{a\}$. Then $d_B(f^0,V,0)\neq0$.

\vskip 0.2 true cm
\noindent
{\bf Lemma 2.2.}\, {\it [6]. Consider $p+1$ linearly independent analytical functions $f_i:U\rightarrow \mathbb{R},\ i=1,2,\cdots, p$, where $U\in\mathbb{R}$ is an interval. Suppose that there exists $j\in\{0,1,\cdots,p\}$ such that $f_j$ has constant sign. Then there exists $p+1$ constants $\delta_i,\ i=0,1,\cdots,p$,  such that $f(x)=\sum\limits_{i=0}^p\delta_if_i(x)$ has at least $p$ simple zeros in $U$. }

\vskip 0.2 true cm

\section{The expression of averaged function and the lower bounds}

Let $x=r\cos\theta$, $y=r\sin\theta$. For $x>0$, system (1.6) can be written as
\begin{eqnarray*}\begin{cases}
\dot{r}=\varepsilon\big(f^+(r,\theta)\cos\theta+g^+(r,\theta)\sin\theta\big),\\
\dot{\theta}=h_1(r,\theta)+\dfrac{\varepsilon}{r}\big(g^+(r,\theta)\cos\theta-f^+(r,\theta)\sin\theta\big),\\
\end{cases}\end{eqnarray*}
where $h_1(r,\theta)=(r \cos\theta +a)^2$. Then, we obtain
$$\frac{dr}{d\theta}=\varepsilon X^+(r,\theta)+\varepsilon^2 Y^+(r,\theta,\varepsilon),\, \cos\theta >0;\eqno(3.1)$$
where
\begin{eqnarray*}\begin{aligned}
X^+(r,\theta)&=[f^+(r,\theta)\cos\theta+g^+(r,\theta)\sin\theta]/h_1(r,\theta),\\[0.2cm]
Y^+(r,\theta,\varepsilon)&=-\dfrac{[f^+(r,\theta)\cos\theta+g^+(r,\theta)\sin\theta][g^+(r,\theta)\cos\theta
-f^+(r,\theta)\sin\theta]}{h_1(r,\theta)[rh_1(r,\theta)+\varepsilon(g^+(r,\theta)\cos\theta-f^+(r,\theta)\sin\theta)]}.
\end{aligned}\end{eqnarray*}
Similarly, we have
$$\frac{dr}{d\theta}=\varepsilon X^-(r,\theta)+\varepsilon^2Y^-(r,\theta,\varepsilon),\, \cos\theta <0,\eqno(3.2)$$
where
\begin{eqnarray*}\begin{aligned}
X^-(r,\theta)&=[f^-(r,\theta)\cos\theta+g^-(r,\theta)\sin\theta]/{h_2(r,\theta)},\\[0.2cm]
Y^-(r,\theta,\varepsilon)&=-\frac{[f^-(r,\theta)\cos\theta+g^-(r,\theta)\sin\theta][g^-(r,\theta)\cos\theta
-f^-(r,\theta)\sin\theta]}{h_2(r,\theta)[rh_2(r,\theta)+\varepsilon(g^-(r,\theta)\cos\theta-f^-(r,\theta)\sin\theta)]},
\end{aligned}\end{eqnarray*}
and $h_2(r,\theta)=(r \cos\theta +b)^2$. It is known that $X^{\pm}$ are well defined in $(0, r_0)$ (where $r_0$ is defined by (1.8)). Hence, system (1.6) is equivalent to
$$\frac{dr}{d\theta}=\varepsilon F^0(\theta,r)+\varepsilon^2 R^0(\theta,r,\varepsilon),$$
with
$$F^0(\theta,r)=F_1(\theta,r)+\textup{sign}(\cos\theta)F_2(\theta,r),$$
$$\quad R^0(\theta,r,\varepsilon)=R_1(\theta,r,\varepsilon)+\textup{sign}(\cos\theta)R_2(\theta,r,\varepsilon),$$
where
$$F_1(r,\theta)=\frac{1}{2}\big[X^+(r,\theta)+X^-(r,\theta)\big],\,$$
$$F_2(r,\theta)=\frac{1}{2}\big[X^+(r,\theta)-X^-(r,\theta)\big],$$
$$R_1(\theta,r,\varepsilon)=\frac{1}{2}\big[Y^+(\theta,r,\varepsilon)+
Y^-(\theta,r,\varepsilon)\big],\, $$
$$R_2(\theta,r,\varepsilon)=\frac{1}{2}\big[Y^+(\theta,r,\varepsilon)-
Y^-(\theta,r,\varepsilon)\big].$$
It is obvious that $F_i(\theta,r)$, $R_i(\theta,r,\varepsilon)$, $i=1,2$ and $\mu(\theta,r)=\cos\theta$ are locally Lipschitz with respect $r$. Since $\mathcal{M}=\{(r,\theta)|r\in(0,r_0),\theta=\pi/2,3\pi/2\}$, the function $\partial \mu(\theta,r)/\partial \theta=-\sin\theta\neq0$ when $(r,\theta)\in\mathcal{M}$. According to Lemma 2.1 and (2.2), we  need only to compute the number of simple zeros of the averaged function
\begin{eqnarray*}\begin{aligned}
&\qquad f^0(r)=\int_0^{2\pi}F^0(\theta,r)d\theta\\
&=\int_{-\frac{\pi}{2}}^\frac{\pi}{2}X^+(\theta,r)d\theta
+\int_{\frac{\pi}{2}}^\frac{3\pi}{2}X^-(\theta,r)d\theta\\
&=\sum\limits_{i+j=1}^{n+1}\sigma_{i,j}r^{i+j-1}\int_{-\frac{\pi}{2}}^{\frac{\pi}{2}}
\frac{\cos^i\theta\sin^{j}\theta}{(r\cos\theta+a)^2}d\theta
+\sum\limits_{i+j=1}^{n+1}\tau_{i,j}r^{i+j-1}\int_{\frac{\pi}{2}}^{\frac{3\pi}{2}}
\frac{\cos^i\theta\sin^{j}\theta}{(r\cos\theta+b)^2}d\theta,
\end{aligned}\quad (3.3)\end{eqnarray*}
where
$$\sigma_{i,j}=a^+_{i-1,j}+b^+_{i,j-1},\ \ \tau_{i,j}=a^-_{i-1,j}+b^-_{i,j-1},\eqno(3.4)$$
here we assume that $a^+_{-1,j}=b^+_{i,-1}=a^-_{-1,j}=b^-_{i,-1}=0$.
Since the coefficients $a^\pm_{i,j}$ and $b^\pm_{i,j}$ are arbitrary for $i+j\geq0$, $\sigma_{i,j}$ and $\tau_{i,j}$ are also arbitrary. Also we note that  the zeros of the function $f^0(r)$ coincide with the zeros of ${\cal F}(r)=rf^0(r)$ in $(0,r_0)$. Hence,
we consider the function ${\cal F}(r)$ instead of $f^0(r)$ in order to simplify further computations.
For $i,j,k\in N \cup \{0\}$, let
 $$ m(k)=\int_{-\frac{\pi}{2}}^\frac{\pi}{2}\cos^k\theta d\theta,\eqno(3.5)$$
 \begin{eqnarray*}
 \,\,\,\qquad \qquad \begin{aligned}
 I_{i,j}(r)&=\int_{-\frac{\pi}{2}}^\frac{\pi}{2}\frac{\cos^i\theta\sin^j\theta}{r\cos\theta+a}d\theta,\ \quad\,
J_{i,j}(r)=\int_{\frac{\pi}{2}}^\frac{3\pi}{2}\frac{\cos^i\theta\sin^j\theta}{r\cos\theta+b}d\theta,\\
A_{i,j}(r)&=\int_{-\frac{\pi}{2}}^\frac{\pi}{2}\frac{\cos^i\theta\sin^j\theta}{(r\cos\theta+a)^2}d\theta,\ \,B_{i,j}(r)=\int_{\frac{\pi}{2}}^\frac{3\pi}{2}\frac{\cos^i\theta\sin^j\theta}{(r\cos\theta+b)^2}d\theta.
\end{aligned}\qquad (3.6)\end{eqnarray*}

\noindent
{\bf Lemma 3.1.}\, {\it The following equalities hold:}
$$rA_{i+1,j}(r)=I_{i,j}(r)-aA_{i,j}(r),\,rB_{i+1,j}(r)=J_{i,j}(r)-bB_{i,j}(r);\eqno(3.7)$$
$$\left\{\begin{array}{l}
A_{i,2j+1}(r)=B_{i,2j+1}(r)=0, \\
A_{i,2l}(r)=\sum\limits_{k=0}^l(-1)^kC_l^kA_{i+2k,0}(r), \\
B_{i,2l}(r)=\sum\limits_{k=0}^l(-1)^kC_l^kB_{i+2k,0}(r);\\
\end{array}\right.\eqno(3.8)$$
$$\left\{\begin{array}{l}
r^iA_{i,0}(r)=(-a)^iA_{0,0}(r)+i(-a)^{i-1}I_{0,0}(r)\\
\qquad \qquad\quad+\sum\limits_{k=0}^{i-2}(k+1)(-a)^k
m(i-k-2)r^{i-k-2},\\
r^iB_{i,0}(r)=(-b)^iB_{0,0}(r)+i(-b)^{i-1}J_{0,0}(r)\\
\qquad \qquad\quad+(-1)^i\sum\limits_{k=0}^{i-2}(k+1)b^km(i-k-2)r^{i-k-2},
\end{array}\right.\eqno(3.9)$$
{\it where  $i,j,k\in N \cup \{0\}$.}

\vskip 0.2 true cm
\noindent
{\bf Proof.}~The formula (3.7) follows easily from the direct computations. Noting that the integrands of $A_{i,j}$ and $B_{i,j}$ in (3.6) are odd functions of $\theta$ and the formula $\sin^{2i}\theta=\sum\limits_{s=0}^iC_i^s(-1)^s\cos^{2s}\theta$, we can get (3.8).

Now we begin to prove the first equality in (3.9) by induction on $i$. For $i,j,k\in N\cup\{0\}$, it follows from [19] that
\begin{eqnarray*}
~~~~~~~~~~~~~~~\begin{cases}
rI_{i,0}(r)=m(i-1)-aI_{i-1,0}(r),\\
rJ_{i,0}(r)=(-1)^{i-1}m(i-1)-bJ_{i-1,0}(r),\\
r^iI_{i,0}(r)=(-a)^iI_{0,0}(r)+\sum\limits_{k=0}^{i-1}(-a)^{i-k-1}r^km(k),\\
r^iJ_{i,0}(r)=(-b)^iJ_{0,0}(r)+(-1)^{i-1}\sum\limits_{k=0}^{i-1}b^{i-k-1}r^km(k).\\
\end{cases}~~~~~~~~~(3.10)
\end{eqnarray*}
The Eqn.(3.8) and the first equality in Eqn.(3.10) imply that the first equality in Eqn.(3.9) holds for $i=1,2$. Suppose that it holds for $i=l$. Then for $i=l+1$, noting that Eqn.(3.10) we have
\begin{eqnarray*}
\begin{aligned}
r^{l+1}A_{l+1,0}(r)=&r^l[I_{l,0}(r)-aA_{l,0}(r)]\\
=&(-a)^lI_{0,0}(r)+\sum\limits_{k=0}^{l-1}(-a)^{l-k-1}r^km(k)+(-a)^{l+1}A_{0,0}(r)\\
&+l(-a)^lI_{0,0}(r)+\sum\limits_{k=0}^{l-2}(k+1)(-a)^{k+1}m(l-k-2)r^{l-k-2}\\
=&(-a)^{l+1}A_{0,0}(r)+(l+1)(-a)^lI_{0,0}(r)+\sum\limits_{k=0}^{l-1}(k+1)(-a)^km(l-k-1)r^{l-k-1},
\end{aligned}
\end{eqnarray*}
which implies that the first equality in Eqn.(3.9) holds for $i=l+1$. The second formula in (3.9) can be proved similarly. This ends the proof.\quad $\lozenge$

\vskip 0.2 true cm

By (3.8), we have
\begin{eqnarray*}
\begin{aligned}
{\cal F}(r)=&\sum\limits_{i+j=1}^{n+1}r^{i+j}\sigma_{i,j}A_{i,j}(r)
 +\sum\limits_{i+j=1}^{n+1}r^{i+j}\tau_{i,j}B_{i,j}(r)\\
=&\sum\limits_{i=1}^{n+1}r^{i}\sum\limits_{j=0}^i\sigma_{i-j,j}A_{i-j,j}(r)+
\sum\limits_{i=1}^{n+1}r^i\sum\limits_{j=0}^i\tau_{i-j,j}B_{i-j,j}(r)\\
=&\sum\limits_{i=1}^{n+1}r^i\sum\limits_{j=0}^{[\frac{i}{2}]}\sigma_{i-2j,2j}A_{i-2j,2j}(r)+
\sum\limits_{i=1}^{n+1}r^i\sum\limits_{j=0}^{[\frac{i}{2}]}\tau_{i-2j,2j}B_{i-2j,2j}(r)\\
=&\sum\limits_{i=0}^{n+1}\sum\limits_{j=0}^{[\frac{n+1-i}{2}]}S_{i,j}r^{i+2j}A_{i,0}(r)+
\sum\limits_{i=0}^{n+1}\sum\limits_{j=0}^{[\frac{n+1-i}{2}]}T_{i,j}r^{i+2j}B_{i,0}(r),
\end{aligned}\end{eqnarray*}
where
$$S_{i,j}=\sum\limits_{k=0}^{[\frac{i}{2}]}(-1)^kC^k_{j+k}\sigma_{i-2k,2j+2k},\ \
T_{i,j}=\sum\limits_{k=0}^{[\frac{i}{2}]}(-1)^kC^k_{j+k}\tau_{i-2k,2j+2k}.\eqno(3.11)$$
By (3.9), we have
\begin{eqnarray*}
\begin{aligned}
{\cal F}(r)=&\sum\limits_{i=0}^{n+1}\sum\limits_{j=0}^{[\frac{n+1-i}{2}]}S_{i,j}r^{2j}[(-a)^iA_{0,0}(r)+i(-a)^{i-1}I_{0,0}(r)\\
&+\sum\limits_{k=0}^{i-2}(k+1)(-a)^km(i-k-2)r^{i-k-2}]\\
&+\sum\limits_{i=0}^{n+1}\sum\limits_{j=0}^{[\frac{n+1-i}{2}]}T_{i,j}r^{2j}[(-b)^iB_{0,0}(r)+i(-b)^{i-1}J_{0,0}(r)\\
&+(-1)^i\sum\limits_{k=0}^{i-2}(k+1)b^km(i-k-2)r^{i-k-2}].
\end{aligned}\end{eqnarray*}
Hence,
\begin{eqnarray*}
\begin{aligned}
{\cal F}(r)=&\sum\limits_{i=0}^{[\frac{n+1}{2}]}\sum\limits_{j=0}^{n+1-2i}S_{j,i}(-a)^jr^{2i}A_{0,0}(r)+
\sum\limits_{i=0}^{[\frac{n+1}{2}]}\sum\limits_{j=1}^{n+1-2i}jS_{j,i}(-a)^{j-1}r^{2i}I_{0,0}(r)\\
&+\sum\limits_{i=0}^{n-1}\sum\limits_{j=0}^{[\frac{i}{2}]}\sum\limits_{k=0}^{n-1-i}S_{i+k-2j+2,j}(k+1)(-a)^km(i-2j)r^i\\
\end{aligned}\end{eqnarray*}
\begin{eqnarray*}
\begin{aligned}
&+\sum\limits_{i=0}^{[\frac{n+1}{2}]}\sum\limits_{j=0}^{n+1-2i}T_{j,i}(-b)^jr^{2i}B_{0,0}(r)+
\sum\limits_{i=0}^{[\frac{n+1}{2}]}\sum\limits_{j=1}^{n+1-2i}jT_{j,i}(-b)^{j-1}r^{2i}J_{0,0}(r)\\
&+\sum\limits_{i=0}^{n-1}\sum\limits_{j=0}^{[\frac{i}{2}]}\sum\limits_{k=0}^{n-1-i}(-1)^iT_{i+k-2j+2,j}(k+1)b^km(i-2j)r^i.
\end{aligned}\end{eqnarray*}
Using the formula
$$\int\frac{1}{(r\cos\theta+a)^2}d\theta=\frac{r\sin\theta}{(r^2-a^2)(r\cos\theta+a)}
-\frac{a}{r^2-a^2}\int\frac{1}{r\cos\theta+a}d\theta,$$
we have
\begin{eqnarray*}\begin{aligned}
&I_{0,0}(r)=\frac{2}{a^2}r+aA_{0,0}(r)-\frac{1}{a}r^2A_{0,0}(r),\\
&J_{0,0}(r)=-\frac{2}{b^2}r+bB_{0,0}(r)-\frac{1}{b}r^2B_{0,0}(r).\end{aligned}\end{eqnarray*}
Thus, we have proved

\vskip 0.2 true cm
\noindent
{\bf Lemma 3.2.} ${\cal F}(r)$ can be expressed as
\begin{eqnarray*}
\begin{aligned}
{\cal F}(r)=\sum\limits_{i=0}^{[\frac{n+1}{2}]+1}a_ir^{2i}A_{0,0}(r)+\sum\limits_{i=0}^{2[\frac{n+1}{2}]+1}b_ir^i
+\sum\limits_{i=0}^{[\frac{n+1}{2}]+1}c_ir^{2i}B_{0,0}(r)+\sum\limits_{i=0}^{2[\frac{n+1}{2}]+1}d_ir^i,
\end{aligned} ~~~~(3.12)
\end{eqnarray*}
where
$$\begin{cases}
a_0=\sum\limits_{j=1}^{n+1}(-a)^j S_{j,0}-\sum\limits_{j=1}^{n+1}j(-a)^jS_{j,0}=-\sum\limits_{j=0}^{n-1}(j+1)(-a)^{j+2}S_{j+2,0},\\[0.4true cm]
a_i=\sum\limits_{j=0}^{n+1-2i}(-a)^j S_{j,i}-\sum\limits_{j=1}^{n+1-2i}j(-a)^jS_{j,i}+
\sum\limits_{j=1}^{n+3-2i}j(-a)^{j-2}S_{j,i-1},\ 1\leq i\leq[\frac{n+1}{2}],\\
a_{[\frac{n+1}{2}]+1}=-\sum\limits_{j=1}^{n+1-2[\frac{n+1}{2}]}j(-a)^{j-2}S_{j,[\frac{n+1}{2}]};
\end{cases}$$
$$\begin{cases}
b_0=\sum\limits_{k=0}^{n-1}(k+1)(-a)^kS_{k+2,0}m(0),\\[0.4true cm]
b_i=\begin{cases}
\sum\limits_{j=0}^{[\frac{i}{2}]}\sum\limits_{k=0}^{n-1-i}(k+1)(-a)^km(i-2j)S_{i+k-2j+2,j}
+\sum\limits_{k=1}^{n+2-i}2k(-a)^{k-3}S_{k,[\frac{i-1}{2}]},\  \textup{if}\ n\ \textup{is \ odd},\\
\sum\limits_{j=0}^{[\frac{i}{2}]}\sum\limits_{k=0}^{n-1-i}(k+1)(-a)^km(i-2j)S_{i+k-2j+2,j},
\qquad\qquad\qquad\qquad\qquad\quad\ \textup{if}\ n\ \textup{is \ even},\\
\end{cases}\\
    \qquad \qquad \qquad \qquad \qquad \qquad\qquad \qquad \qquad 1\leq i\leq2[\frac{n+1}{2}]-1,\\
b_{2[\frac{n+1}{2}]}=0,\ b_{2[\frac{n+1}{2}]+1}=\sum\limits_{j=1}^{n+1-2[\frac{n+1}{2}]}2j(-a)^{j-3}S_{j,[\frac{n+1}{2}]};
\end{cases}$$
$$\begin{cases}
c_0=\sum\limits_{j=1}^{n+1}(-b)^j T_{j,0}-\sum\limits_{j=1}^{n+1}j(-b)^jT_{j,0}=-\sum\limits_{j=0}^{n-1}(j+1)(-b)^{j+2}T_{j+2,0},\\
c_i=\sum\limits_{j=0}^{n+1-2i}(-b)^j T_{j,i}-\sum\limits_{j=1}^{n+1-2i}j(-b)^jT_{j,i}+
\sum\limits_{j=1}^{n+3-2i}j(-b)^{j-2}T_{j,i-1},\ 1\leq i\leq[\frac{n+1}{2}],\\
c_{[\frac{n+1}{2}]+1}=-\sum\limits_{j=1}^{n+1-2[\frac{n+1}{2}]}j(-b)^{j-2}T_{j,[\frac{n+1}{2}]};\\
\end{cases}$$ $$\begin{cases}
d_0=\sum\limits_{k=0}^{n-1}(k+1)(-b)^kT_{k+2,0}m(0),\\
d_i=\begin{cases}
\sum\limits_{j=0}^{[\frac{i}{2}]}\sum\limits_{k=0}^{n-1-i}(-1)^i(k+1)(-b)^km(i-2i)T_{i+k-2j+2,j}
-\sum\limits_{k=1}^{n+2-i}2k(-b)^{k-3}T_{k,[\frac{i-1}{2}]},\  \textup{if}\ n\ \textup{is \ odd},\\
\sum\limits_{j=0}^{[\frac{i}{2}]}\sum\limits_{k=0}^{n-1-i}(-1)^i(k+1)(-b)^km(i-2j)T_{i+k-2j+2,j},
\qquad\qquad\qquad\qquad\qquad\quad\ \textup{if}\ n\ \textup{is \ even},\\
\end{cases}\\
 \qquad \qquad \qquad \qquad \qquad \qquad\qquad \qquad \qquad
1\leq i\leq2[\frac{n+1}{2}]-1,\\
d_{2[\frac{n+1}{2}]}=0,\ d_{2[\frac{n+1}{2}]+1}=\sum\limits_{j=1}^{n+1-2[\frac{n+1}{2}]}2j(-b)^{j-3}T_{j,[\frac{n+1}{2}]}.
\end{cases}$$

\vskip 0.2 true cm
\noindent
{\bf Remark 3.1.} (i)~It is easy to see that
$$a_0=-\frac{a^2}{\pi}b_0,~~ c_0=-\frac{b^2}{\pi}d_0.\eqno(3.13)$$

(ii)~Since $\sigma_{i,j}$ and $\tau_{i,j}$ are arbitrary, (3.11) implies that $S_{i,j}$ and $T_{i,j}$ are also arbitrary.

\vskip 0.2 true cm
\noindent
{\bf Lemma 3.3.} {\it (i)~If $n$ is odd, then the coefficients $a_i$, $b_k$, $c_i$ and $d_k$ in (3.12) are arbitrary for $i=0,1,\cdots,[\frac{n+1}{2}]$ and $k=1,2,\cdots, 2[\frac{n+1}{2}]-1$.

(ii)~If $n$ is even, then the coefficients $a_i$, $b_k$, $c_i$ and $d_k$ in (3.12) are arbitrary for $i=0,1,\cdots,[\frac{n+1}{2}]+1$ and $k=1,2,\cdots,2[\frac{n+1}{2}]-1,2[\frac{n+1}{2}]+1$.}

\vskip 0.2 true cm
\noindent
{\bf Proof.} Suppose that $n$ is odd.  We have
\begin{eqnarray*}
\begin{aligned}
&A=\frac{\partial(a_1,\cdots,a_{[\frac{n+1}{2}]},a_0,b_1,b_2,b_3,b_4,\cdots,b_{2[\frac{n+1}{2}]-2},b_{2[\frac{n+1}{2}]-1})}
{\partial(S_{0,1},\cdots,S_{0,[\frac{n+1}{2}]},S_{1,0},S_{2,0},S_{3,0},S_{5,0},S_{6,0},\cdots,S_{2[\frac{n+1}{2}],0},S_{1,[\frac{n+1}{2}]-1})}\\\\
&=\left(\begin{matrix}
1&\cdots&0&\frac{1}{a}&-2&3a&5a^3&\cdots&-2[\frac{n+1}{2}]a^{2[\frac{n+1}{2}]-2}&0\\
\vdots&\vdots&\vdots&\vdots&\vdots&\vdots&\vdots&\vdots&\vdots&\vdots\\
0&\cdots&1&0&0&0&0&\cdots&0&\frac{1}{a}\\
0&\cdots&0&0&-a^2&2a^3&4a^5&\cdots&-\big(2[\frac{n+1}{2}]-1\big)a^{2[\frac{n+1}{2}]}&0\\
0&\cdots&0&\frac{2}{a^2}&-\frac{4}{a}&\lambda_1&\lambda_2&\cdots&\lambda_3\\
0&\cdots&0&\frac{2}{a^2}&-\frac{4}{a}&6&\lambda_4&\cdots&\lambda_5&0\\
0&\cdots&0&0&0&0&m(3)&\cdots&\lambda_6&0\\
\vdots&\vdots&\vdots&\vdots&\vdots&\vdots&\vdots&\vdots&\vdots&\vdots\\
0&\cdots&0&0&0&0&0&\cdots&m(2[\frac{n+1}{2}]-2)&0\\
0&\cdots&0&0&0&0&0&\cdots&0&\frac{2}{a^2}
\end{matrix}\right),\ \
\end{aligned}
\end{eqnarray*}
where
\begin{eqnarray*}
\begin{aligned}
\lambda_1&=6+m(1),\ \  \lambda_2=10a^2+3a^2m(1), \\
\lambda_3&=-\big(2[\frac{n+1}{2}]-2\big)m(1)a^{2[\frac{n+1}{2}]-3},\\
\lambda_4&=10a^2-2am(2),\\
\lambda_5&=\big(2[\frac{n+1}{2}]-3\big)m(2)a^{2[\frac{n+1}{2}]-4},\\
\lambda_6&=-\big(2[\frac{n+1}{2}]-4\big)m(3)a^{2[\frac{n+1}{2}]-5},
\end{aligned}
\end{eqnarray*}
which implies that $\det A=-\frac{8}{a^2}\prod\limits_{k=3}^{2[\frac{n+1}{2}]-2}m(k)\neq0$. Since $S_{i,j}$ are arbitrary by Remark 3.1(ii), $a_i\ (i=0,1,\cdots,[\frac{n+1}{2}])$ and $b_k\ (k=1,2,\cdots,2[\frac{n+1}{2}]-1)$ can be chosen arbitrarily.
In a similar way, we can prove that $c_i\ (i=0,1,\cdots,[\frac{n+1}{2}])$, $d_k\ (k=1,2,\cdots,2[\frac{n+1}{2}]-1)$ are also arbitrary.

Following the same argument, we can prove the conclusion for $n$ even. This ends the proof. \quad $\lozenge$

\vskip 0.2 true cm

In what follows, we intend to obtain the lower bounds of the number of zeros of ${\cal F}(r)$ by Lemma 2.2.  To this end, we will suppose that $r$ and ${\cal F}(r)$ are complex from now on.

\vskip 0.2 true cm
\noindent
{\bf Lemma 3.4.}\, (i) {\it If $a>0$, then $A_{0,0}(r)$ can be analytically extended to the complex domain $D_1=\mathbb{C}\setminus \{r\in\mathbb{R}|r\leq-a\}$.}

(ii) {\it If $b>0$, then $B_{0,0}(r)$ can be analytically extended to the complex domain $D_2=\mathbb{C}\setminus \{r\in\mathbb{R}|r\geq b\}$.}

\vskip 0.2 true cm
\noindent
{\bf Proof.}~~If $a>0$, then for $\theta\in[-\frac{\pi}{2},\frac{\pi}{2}]$, $\frac{1}{(r\cos\theta+a)^2}$ is analytic for $r\in(-a,+\infty)$. Of course, $A_{0,0}(r)=\int_{-\frac{\pi}{2}}^{\frac{\pi}{2}}\frac{1}{(r\cos\theta+a)^2}d\theta$ is  analytic for $r\in(-a,+\infty)$. By directly computations, we obtain
\begin{eqnarray*}
\begin{aligned}
A_{0,0}(r)=\begin{cases}
\frac{2r}{a(r^2-a^2)}-\frac{4a}{(a^2-r^2)\sqrt{a^2-r^2}}\arctan\sqrt{\frac{a-r}{a+r}},\  a<0,\ r\in(0,-a)\cup(a,0),\\
\infty,\qquad\qquad\qquad\qquad\qquad\qquad\qquad   \quad                                   a<0,\ r=-a,\\
\frac{2r}{a(r^2-a^2)}-\frac{2a}{(r^2-a^2)\sqrt{r^2-a^2}}\ln\frac{r+\sqrt{r^2-a^2}}{-a},\,  \quad a<0,\ r\in(-a,+\infty)\cup(-\infty,a),\\
\frac{2r}{a(r^2-a^2)}+\frac{4a}{(a^2-r^2)\sqrt{a^2-r^2}}\arctan\sqrt{\frac{a-r}{a+r}},\  ~a>0,\ r\in(0,a)\cup(-a,0),\\
\frac{4}{3a^2},\qquad\qquad\qquad\qquad\qquad\qquad\qquad \quad a>0,\ r=a,\\
\frac{2r}{a(r^2-a^2)}-\frac{2a}{(r^2-a^2)\sqrt{r^2-a^2}}\ln\frac{r+\sqrt{r^2-a^2}}{a},\quad\ \, a>0.\ r\in(a,+\infty)\cup(-\infty,-a).\\
\end{cases}
\end{aligned} (3.14)
\end{eqnarray*}
So $A_{0,0}(r)$ satisfies the equation
$$a(a^2-r^2)A'_{0,0}(r)=3arA_{0,0}(r)-4.\eqno(3.15)$$
Solving the above equation, we get
$$A_{0,0}(r)=\frac{1}{(a^2-r^2)\sqrt{a^2-r^2}}
\Big(\pi-2a\int_0^r\frac{1}{\sqrt{a^2-z^2}}dz\Big)-\frac{2r}{a(a^2-r^2)}.$$
Noting that $\frac{1}{\sqrt{a^2-r^2}}$ and $\frac{1}{a^2-r^2}$ is analytic in the domain $\mathbb{C}\setminus\{r\in\mathbb{R}|r^2\geq a^2\}$, we obtain that $A_{0,0}(r)$ is analytic in the domain $\mathbb{C}\setminus\{r\in\mathbb{R}|r^2\geq a^2\}$. Hence, $A_{0,0}(r)$ is analytic in the domain $\mathbb{C}\setminus\{r\in\mathbb{R}|r\leq-a\}$.
The result about $B_{00}$ can be proved similarly. \quad $\lozenge$

\vskip 0.2 true cm

For $r<-a$, denote $A^\pm_{0,0}(r)$ by the analytic continuation of $A_{0,0}(r)$ along an arc such that Im$(r)>0$ (resp.\,Im$(r)<0$). For $r>b$, denote $B^\pm_{0,0}(r)$ by the analytic continuation of $B_{0,0}(r)$ along an arc such that Im$(r)>0$ (resp.\,Im$(r)<0$). For other functions, we will use similar notations.

\vskip 0.2 true cm
\noindent
{\bf Lemma 3.5.}\, {\it If $a>0$, then  $A_{0,0}(r)\thicksim\frac{\pi}{\sqrt{2a}(r+a)^\frac{3}{2}}$ when $r\rightarrow-a$, and
$A_{0,0}(r)\thicksim\frac{2}{ar}$ when $r\rightarrow\infty$. If $b>0$, then  $B_{0,0}(r)\thicksim\frac{\pi}{\sqrt{2b}(r-b)^\frac{3}{2}}$ when $r\rightarrow b$, and
 $B_{0,0}(r)\thicksim\frac{2}{br}$ when $r\rightarrow\infty$.}

\vskip 0.2 true cm
\noindent
{\bf Proof.}\, For $a>0$, from (3.14) we have that
$A_{0,0}(r)\thicksim\frac{\pi}{\sqrt{2a}(a+r)^\frac{3}{2}}$ when $r\rightarrow-a^+$, and
$A_{0,0}(r)\thicksim\frac{2}{ar}$
when $r\rightarrow+\infty$. It follows from (3.15) that
$$(a^2-r^2)A''_{0,0}(r)-5rA'_{0,0}(r)-3A_{0,0}(r)=0.\eqno(3.16)$$
It is easy to check that $-a$ and $\infty$ are singularity points of Fuchs type, so the solution $A_{0,0}(r)$ is regular at these two points. Hence,  $A_{0,0}(r)\thicksim\frac{\pi}{\sqrt{2a}(a+r)^\frac{3}{2}}$ when $r\rightarrow-a$, and  $A_{0,0}(r)\thicksim\frac{2}{ar}$ when $r\rightarrow\infty$.
The proof for the function $B_{0,0}(r)$ is similar, thus we omit it.\quad $\lozenge$

\vskip 0.2 true cm
\noindent
{\bf Lemma 3.6.}\, {\it If $a>0$, $r\in(-\infty,-a)$, then the function $A^\pm_{0,0}(r)$ satisfy
$$A^+_{0,0}(r)-A^-_{0,0}(r)=\frac{c_1\mathbf{i}}{(a^2-r^2)^\frac{3}{2}};\eqno(3.17)$$
 {\it If $b>0$, $r\in(b,+\infty)$, then the function $B^\pm_{0,0}(r)$ satisfy
$$B^+_{0,0}(r)-B^-_{0,0}(r)=\frac{c_2\mathbf{i}}{(b^2-r^2)^\frac{3}{2}},\eqno(3.18)$$
where $c_1$ and $c_2$ are non-zero real constant and $\mathbf{i}^2=-1$.}}

\vskip 0.2 true cm
\noindent
{\bf Proof.} From (3.15) and noting that $A^\pm_{0,0}(r)$ are both analytic continuation of $A_{0,0}(r)$, we have
\begin{eqnarray*}a(a^2-r^2)\big(A^\pm_{0,0}(r)\big)'=3arA^\pm_{0,0}(r)-4,\end{eqnarray*}
which implies
$$(a^2-r^2)\big(A^+_{0,0}(r)-A^-_{0,0}(r)\big)'=3r\big(A^+_{0,0}(r)-A^-_{0,0}(r)\big).$$
Solving the above equation we obtain
$$A^+_{0,0}(r)-A^-_{0,0}(r)=\frac{c}{(a^2-r^2)^\frac{3}{2}},$$
where $c\in\mathbb{C}$ is a constant. Since $A^+_{0,0}(r)$ and $A^-_{0,0}(r)$ conjugate each other, $A^+_{0,0}(r)-A^-_{0,0}(r)$ is a pure imaginary number, and we can suppose $c=c_1\mathbf{i}\ (c_1\in\mathbb{R})$. We claim that $c_1$ is nonzero. Otherwise, $A_{0,0}(r)$ is global single-valued, and $A_{0,0}(r)$ will be analytic at $r=-a$ or $-a$ is a pole of $A_{0,0}(r)$, which is a contradiction with the fact $A_{0,0}(r)\thicksim\frac{\pi}{\sqrt{2a}(a+r)^\frac{3}{2}}$ when $r\rightarrow-a$.
The conclusion (ii) can be shown similarly.\quad $\lozenge$

\vskip 0.2 true cm
\noindent
 {\bf Lemma 3.7.}\, {\it The generating functions of ${\cal F}(r)$ are the the following linearly independent functions.}

(i) For $a\neq-b$ and $n=2k+1$:
$${ \bigl \{  r^i\bigr \}}_{1\le i\le 2k+1}, \, A_{0,0}(r)-\frac{\pi}{a^2},\, \bigl \{ r^{2i}A_{0,0}(r)\bigr \}_{1\le i\le k+1},\, B_{0,0}(r)-\frac{\pi}{b^2}, \, \bigl \{r^{2i}B_{0,0}(r)\bigr \}_{1\le i\le k+1}. $$
\vskip 0.2 true cm
(ii) For $a\neq-b$ and $n=2k$\,:
$$ r, r^2,\ \ r^3,\ \cdots,\ r^{2k-1},\ r^{2k+1},$$
$$A_{0,0}(r)-\frac{\pi}{a^2},\,\bigl \{ r^{2i}A_{0,0}(r)\bigr \}_{1\le i\le k+1},\,  B_{0,0}(r)-\frac{\pi}{b^2}, \, \bigl \{r^{2i}B_{0,0}(r)\bigr \}_{1\le i\le k+1},\,$$
\vskip 0.2 true cm
(iii) For $a=-b$ and $n=2k+1$:
$$\bigl \{ r^i\bigr \}_{1\le i\le 2k+1}, \, A_{0,0}(r)-\frac{\pi}{a^2}, \, \bigl \{ r^{2i}A_{0,0}(r)\bigr \}_{1\le i\le k+1}.$$
\vskip 0.2 true cm
(iv) For $a=-b$ and $n=2k$\,:
$$r, r^2,\ \ r^3,\ \cdots,\ r^{2k-1},\ r^{2k+1}, \,
A_{0,0}(r)-\frac{\pi}{a^2}, \, \bigl \{ r^{2i}A_{0,0}(r)\bigr \}_{1\le i\le k+1}.$$

\vskip 0.2 true cm
\noindent
{\bf Proof.} (i) Suppose that
\begin{eqnarray*}
\begin{aligned}
\Phi(r):=&\alpha_0\big(A_{0,0}(r)-\frac{\pi}{a^2}\big)+\sum\limits_{i=1}^{k+1}\alpha_ir^{2i}A_{0,0}(r)+
\sum\limits_{i=1}^{2k+1}\beta_ir^i\\
&+\gamma_0\big(B_{0,0}(r)-\frac{\pi}{a^2}\big)+\sum\limits_{i=1}^{k+1}\gamma_ir^{2i}B_{0,0}(r)\equiv0,\ r\in(0,r_0).
\end{aligned}
\end{eqnarray*}
We need to prove $\alpha_i=\gamma_i=0\ (i=0,1,\cdots,k+1)$ and $\beta_i=0\ (i=1,2,\cdots,2k+1)$.

If $a>0$, $b>0$, then $\Phi(r)$ can be analytically extended to the domain $D=D_1\cap D_2$. When $r\in(-\infty,-a)$, (3.17) implies
$$\Phi^+(r)-\Phi^-(r)=\sum\limits_{i=0}^{k+1}\alpha_ir^{2i}\big(A^+_{0,0}(r)-A^-_{0,0}(r)\big)
=\frac{c_1\mathbf{i}}{(a^2-r^2)^\frac{3}{2}}\sum\limits_{i=0}^{k+1}\alpha_ir^{2i}.$$
From $\Phi(r)=0$, we have $\sum\limits_{i=0}^{k+1}\alpha_ir^{2i}=0$, which implies that $\alpha_i=0\ (i=0,1,\cdots,k+1)$. When $r\in(b,+\infty)$, we obtain that $\gamma_i=0\ (i=0,1,\cdots,k+1)$ similarly. Hence, $\sum\limits_{i=1}^{2k+1}\beta_ir^i=0$, which implies $\beta_i=0\ (i=1,2,\cdots,2k+1)$.

The other cases, such as  $a>0$ and $b<0$, $a<0$ and $b>0$, and $a<0$ and $b<0$ can be proved similarly.  Following the same argument, we can prove (ii)-(iv).  \quad $\lozenge$

\vskip 0.2 true cm
\noindent
{\bf Remark 3.2}~~By Lemma 2.2 and Lemma 3.7 we can obtain the lower bounds of $H(n)$ given by (1.9) in Theorem 1.1.

\section{The upper bounds}

In the following we will prove the upper bounds of $H(n)$ given in Theorem 1.1 by using the Argument Principle. Suppose that $a>0$, $b>0$ and $n$ is odd.

\vskip 0.2 true cm

From now on we will denote by $R$ a positive real number large enough and $\varepsilon$ a positive real number small enough.
Let $C_R$ be the circle centered at the origin and radius $R$ and consider the points $A, A',B,B'\in C_R$ where $A=(x_A,\varepsilon)$, $A'=(x_A,-\varepsilon)$, $B=(x_B,\varepsilon)$ and $B'=(x_B,-\varepsilon)$. Let $C_{R,\varepsilon}$ be the curve obtained by removing the arcs $AA'$ and $BB'$ of the circle $C_R$, and  let $C_\varepsilon^+$ be the arc $CC'$ of the circle with center at $(b,0)$ and radius $\varepsilon$ and let $C_\varepsilon^-$ be the arc $DD'$ of the circle with center at $(-a,0)$ and radius $\varepsilon$, where $C=(b,\varepsilon)$, $C'=(b,-\varepsilon)$, $D=(-a,\varepsilon)$ and $D'=(-a,-\varepsilon)$. The segment jointing $A$ and $C$ (resp. $A'$ and $C'$) is denoted by $L_b^+$ (resp. $L_b^-$), and the segment jointing $B$ and $D$ (resp. $B'$ and $D'$) is denoted by $L_a^+$ (resp. $L_a^-$). We define the closed curve
 $$G:=C_{R,\varepsilon}\cup C^+_\varepsilon\cup L^\pm_b\cup C^-_\varepsilon \cup L^\pm_a$$
 on the complex plane and denote by $\Omega$ its interior. Consider the counterclockwise orientation on $G$, see Fig.\,1.
Let us use the notation $\rho({\cal F}(r))_{\Gamma}$ to indicate the number of the complete turns of the vector ${\cal F}(r)$ around the path $\Gamma$ in the counterclockwise direction.

On $C_\varepsilon^-$, Lemmas 3.2 and 3.5 we know that ${\cal F}(r)\thicksim\frac{1}{(a+r)^\frac{3}{2}}$. here we neglect the coefficient of $\frac{1}{(a+r)^\frac{3}{2}}$ and in this section we will always do like this. So $\rho ({\cal F}(r))_{C_\varepsilon^-}\leq \frac{3}{2}$ and noting that the number of the complete turns are integer. Hence $\rho ({\cal F}(r))_{C_\varepsilon^-}\leq 1$. Similarly, we have $\rho({\cal F}(r))_{C_\varepsilon^+}\leq 1$.

On $L_a^+\cup L_a^-$, from Lemma 3.6 we have that ${\cal F}(r)$ is real if and only if Im${\cal F}(r)$=0, that is,
$${\cal F}^+(r)-{\cal F}^-(r)=\frac{c_1\mathbf{i}}{(a^2-r^2)^\frac{3}{2}}\sum\limits_{i=0}^{[\frac{n+1}{2}]}a_ir^{2i}=0.$$
So $\rho ({\cal F}(r))_{C_{L_a^+\cup L_a^-}}\leq[\frac{n+1}{2}]$. Similarly, we have $\rho ({\cal F}(r))_{C_{L_b^+\cup L_b^-}}\leq[\frac{n+1}{2}]$.

 On $C_{R,\varepsilon}$, from Lemma 3.5 we have $r^{2[\frac{n+1}{2}]}A_{0,0}(r)$, $r^{2[\frac{n+1}{2}]}B_{0,0}(r)\thicksim r^{2[\frac{n+1}{2}]-1}$, which implies that $\rho({\cal F}(r))_{C_{R,\varepsilon}}\leq 2[\frac{n+1}{2}]-1$.

Combining the above cases, we have
$$\rho({\cal F}(r))_{G}\leq 4[\frac{n+1}{2}]+1.$$
By the Argument Principle, we obtain that ${\cal F}(r)$ has at most $4[\frac{n+1}{2}]+1$ zeros in $\Omega$. Noting that ${\cal F}(0)=0$, we have
$$H(n)\leq 4[\frac{n+1}{2}].$$
We can prove results for other cases similarly.  This ends the proof of Theorem 1.1.

\section{Proof of the Theorem 1.2}
 \setcounter{equation}{0}
\renewcommand\theequation{5.\arabic{equation}}

For the sake of clearness, we list the perturbed smooth integrable differential system as following
\begin{eqnarray*}
\begin{cases}
\dot{x}=-y(x+a)^2+\varepsilon f(x,y),\\
\dot{y}=x(x+a)^2+\varepsilon g(x,y),\\
\end{cases}
\end{eqnarray*}
where $a\neq0$, $0<|\varepsilon|\ll1$,
$$  f(x,y)=\sum\limits_{i+j=0}^na_{i,j}x^iy^j,\ \
        g(x,y)=\sum\limits_{i+j=0}^nb_{i,j}x^iy^j.$$
Similar to (3.5) and (3.6), we define
\begin{eqnarray}
\begin{aligned}
&U_{i,j}(r)=\int_0^{2\pi}\frac{\cos^i\theta\sin^j\theta}{r\cos\theta+a}d\theta,\ \
V_{i,j}(r)=\int_0^{2\pi}\frac{\cos^i\theta\sin^j\theta}{(r\cos\theta+a)^2}d\theta,\\
&p(k)=\int_0^{2\pi}\cos^k\theta d\theta.
\end{aligned}
\end{eqnarray}
It is easy to check that
\begin{eqnarray}
p(k)=\begin{cases}
2\pi,\ \qquad\ ~\textup{if}\ k=0,\\
2\pi\frac{(k-1)!!}{k!!},\ \ \textup{if}\ k\ \textup{is even and }k\geq2\,,\\
0,\ \qquad\quad~~  \textup{otherwise},\\
\end{cases}
\end{eqnarray}
where $k!!=\prod\limits_{i=0}^{[\frac{k}{2}]-1}(k-2i)$, $0!!=1!!=1$.
Using the similar arguments as in the proof of Lemma 3.1, we obtain the following  equalities
\begin{eqnarray}
\begin{cases}
rV_{i+1,j}(r)=U_{i,j}(r)-aV_{i,j}(r),\\
V_{i,2j+1}(r)=0, \\
V_{i,2l}(r)=\sum\limits_{k=0}^l(-1)^kC_l^kV_{i+2k,0}(r), \\
r^iV_{i,0}(r)=(-a)^iV_{0,0}(r)+i(-a)^{i-1}U_{0,0}(r)\\
\qquad \qquad\quad+\sum\limits_{k=0}^{i-2}(k+1)(-a)^k
p(i-k-2)r^{i-k-2},\\
\end{cases}
\end{eqnarray}
where $i,j,k\in N \cup \{0\}$.

From (5.1)-(5.3) and noting that $$U_{0,0}(r)=aV_{0,0}(r)-\frac{r^2}{a}V_{0,0}(r),$$ we have
\begin{eqnarray*}
\begin{aligned}
{\cal F}(r)=\sum\limits_{i=0}^{[\frac{n+1}{2}]+1}\alpha_ir^{2i}V_{0,0}(r)+\sum\limits_{i=0}^{[\frac{n-1}{2}]}\beta_{2i}r^{2i},
\end{aligned}
\end{eqnarray*}
where
$$\begin{cases}
\alpha_0=-\sum\limits_{j=0}^{n-1}(j+1)(-a)^{j+2}Q_{j+2,0},\\
\alpha_i=\sum\limits_{j=0}^{n+1-2i}(-a)^j Q_{j,i}-\sum\limits_{j=1}^{n+1-2i}j(-a)^jQ_{j,i}+
\sum\limits_{j=1}^{n+3-2i}j(-a)^{j-2}Q_{j,i-1},\ 1\leq i\leq[\frac{n+1}{2}],\\
\alpha_{[\frac{n+1}{2}]+1}=-\sum\limits_{j=1}^{n+1-2[\frac{n+1}{2}]}j(-a)^{j-2}Q_{j,[\frac{n+1}{2}]},\\
\beta_0=\sum\limits_{k=0}^{n-1}(k+1)(-a)^kQ_{k+2,0}p(0),\\
\beta_{2i}=\sum\limits_{j=0}^{i}\sum\limits_{k=0}^{n-1-2i}(k+1)(-a)^kp(2i-2j)Q_{2i+k-2j+2,j},
1\leq i\leq[\frac{n-1}{2}],
\end{cases}$$
and
$$
Q_{i,j}=\sum\limits_{k=0}^{[\frac{i}{2}]}(-1)^kC^k_{j+k}\lambda_{i-2k,2j+2k},\ \
\lambda_{i,j}=a_{i-1,j}+b_{i,j-1},\ \ a_{-1,j}=b_{i,-1}=0.$$

Just following the idea in the proof of Lemma 3.7, we can obtain the following lemma.
\vskip 0.2 true cm
\noindent
 {\bf Lemma 5.1.}\, {\it The generating functions of ${\cal F}(r)$ are the the following linearly independent functions.}

(i) For $n=2k+1$:
$${\bigl \{  r^{2i}\bigr \}}_{1\le i\le k}, \,\ V_{0,0}(r)-\frac{2\pi}{a^2},\,\ \bigl \{ r^{2i}V_{0,0}(r)\bigr \}_{1\le i\le k+1}. $$
\vskip 0.2 true cm
(ii) For $n=2k$\,:
$${\bigl \{  r^{2i}\bigr \}}_{1\le i\le k}, \,\ V_{0,0}(r)-\frac{2\pi}{a^2},\,\ \bigl \{ r^{2i}V_{0,0}(r)\bigr \}_{1\le i\le k+1}. $$

Hence, by Lemma 2.2 and Argument Principle, we can prove Theorem 1.2. This ends the proof.

 \vskip 0.2 true cm

\end{document}